\documentclass[letterpaper, 10 pt, conference]{ieeeconf}  

\IEEEoverridecommandlockouts
\let\labelindent\relax

\usepackage{enumitem}
\usepackage{graphicx}
\usepackage{amsthm,amsmath,amssymb,amsfonts}
\usepackage{grffile}
\usepackage{multicol}
\usepackage{algorithm,algorithmic}
\usepackage{empheq}
\usepackage{mathtools}
\usepackage{xcolor}

\setlist[enumerate]{wide=0pt, widest=99,leftmargin=\parindent, labelsep=*}

\newcommand{\mcl}[1]{\mathcal{ #1}}
\newcommand{\mbf}[1]{\mathbf{ #1}}
\newcommand{\norm}[1]{\left\Vert #1\right\Vert}
\newcommand{\hinf}{\ensuremath{H_{\infty}}}
\newcommand{\ip}[2]{\left\langle{#1},{#2}\right\rangle}

\newcommand{\bmat}[1]{\begin{bmatrix} #1\end{bmatrix}}

\newcommand{\mat}[1]{\begin{matrix}#1\end{matrix}}

\newcommand{\E}{\mathbb{E}}
\newcommand{\R}{\mathbb{R}}

\newcommand{\Z}{\mathbb{Z}}

\newcommand{\myint}{\int_{a}^{b}}

\newtheorem{thm}{Theorem}
\newtheorem{defn}[thm]{Definition}
\newtheorem{lem}[thm]{Lemma}

\newtheorem{cor}[thm]{Corollary}

\newtheorem{ex}[thm]{\textbf{Example}}

\newtheorem{pse}{\textbf{Pseudo Code}}

\newcommand{\fourpi}[4]{\ensuremath{\mcl{P}{\tiny\bmat{#1,& \hspace{-3mm}#2 \\ #3,& \hspace{-3mm} \left\{#4\right\} }}}}
\newcommand{\threepi}[1]{\mcl{P}_{\{#1\}}}

\newcommand{\threepiFull}[3]{\mcl{P}_{\{#1,#2,#3\}}}

\allowdisplaybreaks

\title{\LARGE \bf
Duality and $\hinf$-Optimal Control Of Coupled ODE-PDE Systems
}

\author{Sachin Shivakumar$^{1}$ \and Amritam Das$^{2}$ \and Siep Weiland$^{2}$ \and Matthew M. Peet$^{1}$
\thanks{$^{1}$ Sachin Shivakumar\{sshivak8@asu.edu\} and Matthew M. Peet\{mpeet@asu.edu\} are with School for Engineering of Matter, Transport and Energy, Arizona State University, Tempe, AZ, 85298 USA}%
\thanks{$^{2}$ Amritam Das\{am.das@tue.nl\} and Siep Weiland\{S.Weiland@tue.nl\} are with Department of Electrical Engineering, Eindhoven University of Technology}%
}

\begin{document}

\maketitle
\thispagestyle{empty}
\pagestyle{empty}

\begin{abstract}
	
	In this paper, we present a convex formulation of $\hinf$-optimal control problem for coupled linear ODE-PDE systems with one spatial dimension. First, we reformulate the coupled ODE-PDE system as a Partial Integral Equation (PIE) system and show that stability and $\hinf$ performance of the PIE system implies that of the ODE-PDE system. We then construct a dual PIE system and show that asymptotic stability and $\hinf$ performance of the dual system is equivalent to that of the primal PIE system. Next, we pose a convex dual formulation of the stability and $\hinf$-performance problems using the Linear PI Inequality (LPI) framework. LPIs are a generalization of LMIs to Partial Integral (PI) operators and can be solved using PIETOOLS, a MATLAB toolbox. Next, we use our duality results to formulate the stabilization and $\hinf$-optimal state-feedback control problems as LPIs. Finally, we illustrate the accuracy and scalability of the algorithms by constructing controllers for several numerical examples.
\end{abstract}

\section{INTRODUCTION}\label{sec:introduction}
In this paper, we consider the problem of $\hinf$-optimal state-feedback controller synthesis for Partial Integral Equation (PIE) systems of the form
\begin{align}\label{eq:PIE_general}
\mcl T \dot{\mbf x}(t)&=\mcl A\mbf x(t)+\mcl{B}w(t), \quad\mbf x(0)=\mbf x_0\in \R^m \times L_2^n\notag\\
z(t) &= \mcl{C}\mbf x(t) + \mcl{D}w(t)
\end{align}
where $\mcl{T}, \mcl{A}, \mcl{B}, \mcl{C}, \mcl{D}$ are Partial Integral (PI) operators and $w(t)\in \R^p$. The dual (or adjoint) PIE system is then defined to be
\begin{align}\label{eq:PIE_adjoint_general}
\mcl T^* \dot{\bar{\mbf x}}(t)&=\mcl A^*\bar{\mbf x}(t)+\mcl{C}^*\bar{w}(t), \quad \bar{\mbf x}(0)=\bar{\mbf x}_0\in \R^m \times L_2^n\notag\\
\bar{z}(t) &= \mcl{B}^*\bar{\mbf x}(t) + \mcl{D}^*\bar{w}(t)
\end{align}
where $^*$ denotes the adjoint with respect to $L_2$-inner product. Recently, it has been shown that almost any PDE system in a single spatial dimension coupled with an ODE at the boundary has an equivalent PIE system representation \cite{shivakumar2019generalized} (see Sec. \ref{subsec:ODEPDE}). It should be noted, however, that the formulation in Eqn.~\eqref{eq:PIE_general} does not allow for inputs directly at the boundary - rather these must enter through the ODE or into the domain of the PDE.  Use of the PIE system representation, defined by the algebra of Partial Integral (PI) operators, allows us to generalize LMIs developed for ODEs to infinite-dimensional systems. These generalizations are referred to as Linear PI Inequalities (LPIs) and can be solved efficiently using the Matlab toolbox PIETOOLS~\cite{toolbox:pietools}.  In previous work, LPIs have been proposed for stability~\cite{PEET_CDPS}, $\hinf$-gain~\cite{shivakumar2019computing} and $\hinf$-optimal estimation~\cite{das_2019CDC} of PIE systems. However, until now the stabilization and $\hinf$-optimal controller synthesis problems have remained unresolved. In this paper, we resolve the problems of stabilizing and $\hinf$ state-feedback controller synthesis by proving the following results.
\begin{enumerate}[label=(\Alph*),ref=\Alph*]
	\item \label{item:A}~~\textit{Dual Stability Theorem:} We show that the PIE system \eqref{eq:PIE_general} is stable for $w=0$ and any initial condition $\mbf x(0)\in L_2$ if and only if the dual PIE system \eqref{eq:PIE_adjoint_general} is stable for any initial conditions $\mbf{\bar{x}}(0)\in L_2$ and $\bar{w}=0$.
	\item \label{item:B}~~\textit{Dual $L_2$-gain Theorem:} For $w\in L_2([0,\infty))$ and $\mbf{x}(0)=0$,
	any solution of the PIE system~\eqref{eq:PIE_general} satisfies $\norm{z}_{L_2}\le\gamma\norm{w}_{L_2}$ if and only if any solution to the dual PIE system Eq. \eqref{eq:PIE_adjoint_general} satisfies $\norm{\bar{z}}_{L_2}\le \gamma\norm{\bar{w}}_{L_2}$ for $\bar{\mbf x}(0)=0$ and $\bar{w}\in L_2([0,\infty))$.
	\item \label{item:C}~\textit{$\hinf$-optimal Control of PIEs:} The stabilization and $H_{\infty}$-optimal state-feedback controller synthesis problem for PIE systems~\eqref{eq:PIE_general} may be formulated as an LPI.
\end{enumerate}

Previous work on controller synthesis for coupled ODE-PDE systems includes the well-established method of backstepping (See e.g.~\cite{krstic2008boundary}) and reduced basis methods (See e.g.~\cite{ito1998reduced}). In the former case, backstepping methods allow for inputs at the boundary and are guaranteed to find a stabilizing controller if one exists. However, the resulting controllers are not optimal in any sense. In the latter case, $\hinf$-optimal controllers are designed for an ODE approximation of the coupled ODE-PDE system. However, these controllers do not have provable performance properties when applied to the actual ODE-PDE, i.e. the $\hinf$-norm of the ODE-PDE system is not same as the $\hinf$-norm of the ODE approximation and indeed, the resulting closed-loop system is often unstable.

The fundamental issue in controller synthesis for both finite-dimensional and infinite-dimensional systems is one of non-convexity. In simple terms, for either a finite or infinite-dimensional system of the form
\[ \dot{x}(t) = \mcl A x(t) + \mcl B u(t),
\]
finding a stabilizing control $u(t)=\mcl Kx(t)$ and a corresponding Lyapunov functional $V(t) = \ip{x(t)}{\mcl Px(t)}_X$ with negative time-derivative gives rise to a bilinear problem in variables $\mcl K$ and $\mcl P$ of the form $(\mcl A+\mcl{BK})^{*}\mcl P+\mcl P(\mcl A+\mcl{BK})\le 0$.

In case of finite-dimensional linear systems, the linear operators $\mcl{P}, \mcl{K}, \mcl{A}$ and $\mcl{B}$ are matrices $P$, $K$, $A$ and $B$. In absence of a controller, the Lyapunov stability test (referred to as primal stability test) can be written as an LMI in positive matrix variable $P>0$ such that $A^TP+PA\le0$. In finite-dimensions, the eigenvalues of $A$ and $A^*$ are the same and hence there is an equivalent dual Lyapunov inequality of the form $AP+PA^T\le0$. Then the test for existence of a stabilizing controller $K$ and a Lyapunov functional $P$ which proves the stability of the closed-loop system can now be written as: find $P>0$ such that $(A+BK)P+P(A+BK)^T\le0$. The key difference, however, is the bilinearity can now be eliminated by introducing new variable $Z=KP$ which leads to the LMI constraint $AP+BZ+(AP+BZ)^T\le0$.

For infinite-dimensional systems, Theorem $5.1.3$ of \cite{Curtain:1995:IIL:207416} is similar to primal stability test for ODEs. The result is similar in the sense that matrices in the constraints of primal stability test for ODE are replaced by linear operators for infinite-dimensional systems, i.e. a test for existence of a positive operator $\mcl{P}>0$ that satisfies the operator-valued constraint $\mcl{A}^*\mcl{P}+\mcl{PA}\le0$. However, there does not exist a dual form of the primal stability test for infinite-dimensional systems. In \cite{peet_2019TAC}, a dual Lyapunov criterion for stability in infinite-dimensional systems was presented. However, the result was restricted to infinite-dimensional systems of the form \[\dot{\mbf x}(t)  =\mcl{A}\mbf{x}(t)+\mcl{B}u(t)\] and included constraints on the image of the operator $\mcl P$ of the form $\mcl{P}(X) = X$ where $X=D(\mcl A)$ is the domain of the infinitesimal generator $\mcl A$. Furthermore, because $\mcl A$ for PDEs is a differential operator, this approach provides no way of enforcing negativity of the dual stability condition. These difficulties in analysis and controller synthesis for PDE systems led to the development of the PIE formulation of the problem - wherein both system parameters $\mcl A, \mcl B,\mcl C,\mcl D, \mcl T$ and the Lyapunov parameter $\mcl P$ lie in the algebra of bounded linear PI operators.


In this work, we adopt the PIE formulation of the ODE-PDE system and propose dual stability and performance tests wherein all operators lie in the PI algebra and do not include additional constraints such as $\mcl P(X)=X$. Specifically, the results (\ref{item:A}) and (\ref{item:B}) lead to LPIs which, by allowing for the variable change trick used in finite-dimensional systems, allows us to propose convex and testable formulations of the stabilization and optimal control problems - resulting in stabilizing or $H_{\infty}$ optimal controllers for coupled PDE-ODE systems where the inputs enter through the ODE or in the domain. More specifically, these methods apply for linear ODE-PDE systems in $1$ spatial variable with a very general set of boundary conditions including Dirichlet, Neumann, Robin, Sturm-Lioville et c. The resulting LPIs are solved numerically using PIETOOLS \cite{toolbox:pietools}, an open-source MATLAB toolbox to handle PI variables and setup PI operator-valued optimization problems. Finally, we note that this is the first result to achieve $H_\infty$-optimal control of coupled ODE-PDE systems. Although we are currently restricted to inputs using an ODE filter or in-domain, we believe the duality results presented here can ultimately be extended to cover inputs applied directly at the boundary.

The paper is organized as follows. After introducing preliminary notations in Section II, in Section III and IV, we introduce the general form of PIE and ODE-PDE under consideration. In Section V, we define the conditions under which PIE and ODE-PDE as equivalent followed by equivalence in stability and $\hinf$-gain in Section VI. Section VII discusses the properties of adjoint PIE systems. In Section VIII and IX, we derive the dual stability theorem and dual $\hinf$-gain theorem for PIEs. Sections XI through XIV present the LPIs developed using dual stability theorem and dual $\hinf$-gain theorem. Examples are illustrated in Section XV and followed by conclusions in Section XVI.

\section{Notation}\label{sec:notation}
We use the calligraphic font, for example $\mcl{A}$, to represent linear operators on Hilbert spaces and the bold font, $\mbf x$, to denote functions in $L_2^n[a,b]$ the set of all square-integrable functions on the domain $[a,b]\subset \R$. The Sobolev space $W_{2,k}[a,b]$ is defined as \begin{align*}
W_{2,k}[a,b] := \{f\in L_2[a,b]\mid \frac{\partial^n f}{\partial s^n}\in L_2[a,b] ~\text{for all} ~n\le k\}.
\end{align*} $Z^{m,n}[a,b]$ denotes the space $\R^{m}\times L_2^{n}[a,b]$ which is equipped with the inner-product
{\small\begin{align*}\ip{\bmat{x_1\\\mbf x_2}}{\bmat{y_1\\\mbf y_2}}_{Z} = x_1^Ty_1 + \ip{\mbf x_2}{\mbf y_2}_{L_2}, \quad \bmat{x_1\\\mbf x_2}, \bmat{y_1\\\mbf y_2}\in Z^{m,n}.
	\end{align*}}
We use $\mbf x_s$ to denote partial derivative of $\frac{\partial\mbf x}{\partial s}$ where the number of repetitions of the subscript $s$ corresponds to the order of the partial derivative and $\dot{\mbf x}$ to denote the partial derivative $\frac{\partial \mbf x}{\partial t}$.

\section{Partial Integral Equations} \label{sec:PIE}

In this section, we will define a PIE system with inputs and disturbances of the form
\begin{align}\label{eq:PIE_full}
\mcl T \dot{\mbf x}(t)&=\mcl A\mbf x(t)+\mcl{B}_1w(t)+\mcl B_2u(t), \quad \mbf x(0)\in Z^{m,n}[a,b]\notag\\
z(t) &= \mcl{C}\mbf x(t) + \mcl{D}_{11}w(t) + \mcl{D}_{12}u(t),
\end{align}
where the $\mcl T, \mcl A: Z^{m,n}[a,b]\to Z^{m,n}[a,b]$, $\mcl B_1:\R^{q}\to Z^{m,n}[a,b]$, $\mcl B_2:\R^{p}\to Z^{m,n}[a,b]$, $\mcl{C}:Z^{m,n}[a,b]\to \R^{r}$, $\mcl{D}_{11}\in\R^{r\times q}$ and $\mcl{D}_{12}\in\R^{r\times p}$ are Partial Integral (PI) operators, defined as follows.

\begin{defn}\label{def:4PI}
	(PI Operators:)	A 4-PI operator is a bounded linear operator between $Z^{m,n}[a,b]$ and $Z^{p,q}[a,b]$ of the form
	\begin{align}\label{eq:4pi}
	&\fourpi{P}{Q_1}{Q_2}{R_i}\bmat{x\\\mbf{y}}(s) = \bmat{Px + \int_{a}^{b}Q_1(s)\mbf{y}(s)ds\\Q_2(s)x+\threepi{R_i}\mbf{y} (s)}
	\end{align}
	where $P\in\R^{p\times m}$ is a matrix, $Q_1:[a,b]\to\R^{p\times n}$, $Q_2:[a,b]\to\R^{q\times m}$ are bounded integrable functions and $\threepi{R_i}:L_2^n[a,b]\to L_2^q[a,b]$ is a 3-PI operator of the form
	\begin{align*}
	&\left(\mcl{P}_{\{R_i\}}\mbf x\right)(s):= \\
	&R_0(s) \mbf x(s) +\int_{a}^s R_1(s,\theta)\mbf x(\theta)d \theta+\int_s^bR_2(s,\theta)\mbf x(\theta)d \theta.
	\end{align*}
\end{defn}

\begin{defn}
	For given $u\in L_2([0,\infty);\R^p)$, $w\in L_2([0,\infty);\R^q)$ and initial conditions $\mbf x_0\in Z^{m,n}[a,b]$, we say that $\mbf x:[0,\infty)\to Z^{m,n}[a,b]$ and $z:[0,\infty)\to\R^r$ satisfy the PIE \eqref{eq:PIE_general} defined by $\{\mcl{T},\mcl{A},\mcl{B},\mcl{C},\mcl{D}\}$ if $\mbf x$ is Fr\'echet differentiable almost everywhere on $[0,\infty)$, $\mbf x(0) = \mbf x_0$ and the equations \eqref{eq:PIE_full} are satisfied for almost all $t\ge0$.
\end{defn}

\section{A General Class of Linear ODE-PDE Systems}\label{subsec:ODEPDE}
In this paper, we consider control of the following class of coupled linear ODE-PDE systems in a single spatial variable $s\in[a,b]$.

{\small
	\begin{align}\label{eq:pde_general}
	&\bmat{\dot{x}(t)\\\dot{\mbf{x}}(s,t)} = \bmat{Ax(t) + \left(\mcl{E}\mbf{x}\right)(t)\\E(s)x(t)+\left(\mcl{A}_p\mbf{x}\right)(s,t)} + \bmat{B_1\\B_2(s)}u(t), \notag\\
	&z(t) = \left(\mcl{C}\bmat{x\\\mbf x}\right)(t) + D u(t),\notag\\
	&B\bmat{\mbf x_c(a,t)\\\mbf x_c(b,t)} = B_xx(t), \qquad \bmat{x(0)\\\mbf x(\cdot,0)}=\mbf x_0\in D(\mcl{A}_d)\notag\\
	&\mbf x(s,t) = \bmat{\mbf x_1(s,t)\\\mbf x_2(s,t)\\\mbf x_3(s,t)}, ~~\mbf x_c(s,t) = \bmat{\mbf x_2(s,t)\\\mbf x_3(s,t)\\\mbf x_{3s}(s,t)}
	\end{align}
}
where the differential operator $\mcl{A}_p$ and the operators $\mcl{C}, \mcl{E}$ are defined as
{\small
	\begin{align}
	&\left(\mcl{A}_p\mbf x\right)(s,t):=\notag\\
	& A_0(s)\bmat{\mbf x_1(s,t)\\\mbf x_2(s,t)\\\mbf x_3(s,t)}+A_1(s)\bmat{\mbf x_{2s}(s,t)\\\mbf x_{3s}(s,t)}+A_2(s)\mbf x_{3ss}(s,t),\notag\\
	&\left(\mcl{E}\mbf x\right)(t):= E_{10}\bmat{\mbf x_c(a,t)\\\mbf x_c(b,t)}+\int_{a}^{b}E_a(s)\bmat{\mbf x_1(s,t)\\\mbf x_2(s,t)\\\mbf x_3(s,t)}ds\notag\\
	&\hspace{4cm}+\int_{a}^{b}E_b(s)\bmat{\mbf x_{2s}(s,t)\\\mbf x_{3s}(s,t)}ds,\notag\\
	&\left(\mcl{C}\bmat{x\\\mbf x}\right)(t):= Cx(t)+C_{10}\bmat{\mbf x_c(a,t)\\\mbf x_c(b,t)}\notag\\
	&\hspace{1cm}+\int_{a}^{b}C_a(s)\bmat{\mbf x_1(s,t)\\\mbf x_2(s,t)\\\mbf x_3(s,t)}ds+\int_{a}^{b}C_b(s)\bmat{\mbf x_{2s}(s,t)\\\mbf x_{3s}(s,t)}ds,\notag
	\end{align}}
and where
{\small
	\begin{align}\label{eq:domain}
	&D(\mcl{A}_d):= \notag\\
	&\left\{\mat{\bmat{x\\\mbf x_1\\\mbf x_2\\\mbf x_3}\in\R^{n_o}\times L_2^{n_1}[a,b]\times W_{2,1}^{n_2}[a,b]\times W_{2,2}^{n_3}[a,b] : \\ B\bmat{x_c(a)\\x_c(b)}=B_xx,   ~\text{where}~x_c(s) = \bmat{\mbf x_2(s)\\\mbf x_3(s)\\\mbf x_{3s}(s)}}  \right\}
	\end{align}}
The ODE states are $x(t)\in\R^{n_o}$, while the PDE states are $\mbf x_i(s,t)\in \R^{n_i}$. The total number of PDE is thus defined to be $n_p = n_1+n_2+n_3$. The ODE-PDE system is defined by the parameters $A_0:[a,b]\to \R^{n_p\times n_p}$, $A_1:[a,b]\to\R^{n_p\times (n2+n3)}$, $A_2:[a,b]\to\R^{n_p\times n_3}$, $E:[a,b]\to\R^{n_p\times n_o}$, $E_a:[a,b]\to\R^{n_o\times n_p}$, $E_b:[a,b]\to \R^{n_o \times (n_2+n_3)}$, $C_a:[a,b]\to\R^{n_z\times n_p}$, $C_b:[a,b]\to\R^{n_z\times (n_2+n_3)}$ and $B_2:[a,b]\to \R^{n_p\times n_u}$ are bounded integrable functions. $A\in\R^{n_o\times n_o}$, $E_{10}\in\R^{n_o\times 2n_r}$,  $C_{10}\in\R^{n_z\times 2n_r}$,  $B_1\in\R^{n_o\times n_u}$,  $D\in\R^{n_z\times n_u}$, $B\in\R^{n_r\times 2n_r}$ has row rank $n_r:=\text{rank}(B)=n_2+2n_3$ and $B_x\in\R^{n_r\times n_o}$. This class of systems includes almost all coupled linear ODE-PDE systems with the constraint that the input does not directly act at the boundary, but rather through the ODE or in the domain of the PDE. The model can also be extended if higher-order spatial derivatives are required.

\textbf{Illustrative Example} To illustrate how this representation is applied to a typical ODE-PDE model, we consider a wave equation coupled with an ODE as shown below.
\begin{align}
&\dot{x}(t) = ax(t) + dw(1,t),\label{eq:ex1}\\
&\ddot{w}(s,t) = cw_{ss}(s,t),\qquad  w(0,t) = kx(t), w_s(1,t) = 0,\notag
\end{align}
where $w(s,t)$ is transverse displacement of the string and $x$ is the ODE state.
These equations may be rewritten in the form \eqref{eq:pde_general}

\begin{align*}
&\dot{x}(t) = ax(t)+d\mbf x_3(1,t),\\
&\bmat{\dot{\mbf x}_1\\\dot{\mbf x}_3}(s,t) = \bmat{0&0\\1&0}\bmat{\mbf x_1\\\mbf x_3}(s,t)+\bmat{c\\0}\mbf x_{3ss}(s,t)\\
&\mbf x_3(0,t) = kx(t), \mbf x_{3s}(1,t) = 0
\end{align*} where $\mbf x_1 = \dot{w}$ and $\mbf x_3 = w$.
The parameters that define the ODE-PDE \eqref{eq:pde_general} are
\begin{align*}
&A=a, A_0 = \bmat{0~0\\1~0}, A_2 =\bmat{c\\0}, E_{10} = \bmat{0~0~d~0}, \\
&B = \bmat{1~0~0~0\\0~0~0~1}, B_x = \bmat{k\\0},
\end{align*} and the rest of the system parameters are zero.

\begin{defn}
	For given $u\in L_2([0,\infty);\R^{n_u})$ and initial conditions $\mbf x_0\in D(\mcl{A}_d)$ as defined in \eqref{eq:domain}, we say that $x:[0,\infty)\to \R^{n_o}$, $\mbf x:[0,\infty)\to L_2^{n_1}[a,b]\times W_{2,1}^{n_2}[a,b]\times W_{2,2}^{n_3}[a,b]$ and $z:[0,\infty)\to\R^{n_z}$ satisfy the ODE-PDE \eqref{eq:pde_general} defined by  $\{A,A_i,B_i,B,B_x,C_{10},C_a,C_b,D,E,E_{10},E_a,E_b\}$ if $x$ is differentiable and $\mbf x$ is Fr\'echet differentiable almost everywhere on $[0,\infty)$, $\bmat{x\\\mbf x}(0)=\mbf x_0$, $\bmat{x\\\mbf x}(t)\in D(\mcl{A}_d)$ and Equations~\eqref{eq:pde_general} hold for almost all $t\ge0$.
\end{defn}


\section{PIE Representation of the ODE-PDE System}
A coupled ODE-PDE of the form \eqref{eq:pde_general} can be written as a PIE system. Furthermore, the solutions of the PIE define solutions of the ODE-PDE and vice-versa. The conversion formulae are given in the appendix in Eqns.~\eqref{eq:PIE_conversion}.
\begin{thm}\label{thm:pde2pie}
	For given $u\in L_2([0,\infty);\R^{n_u})$ and initial conditions $\bmat{x_0\\\mbf x_0} \in D(\mcl{A}_d)$ as defined in \eqref{eq:domain}, suppose $x$, $\mbf{x}$ and $z$ satisfy the ODE-PDE defined by $\{A,A_i,B_i,B,B_x,C_{10},C_a,C_b,D,E,E_{10},E_a,E_b\}$. Then $z$ also satisfies the PIE
	\begin{align*}
	\mcl T \dot{\mbf v}(t)&=\mcl A\mbf v(t)+\mcl{B}u(t), \quad \mbf v(0)=\mbf v_0\notag\\
	z(t) &= \mcl{C}\mbf v(t) + \mcl{D}u(t),
	\end{align*}
	with \begin{align*}
		\mbf v_0 = \bmat{x_0\\\mbf x_{01}\\\mbf{x}_{02,s}\\\mbf{x}_{03,ss}}, \quad \mbf v(t) := \bmat{x\\\mbf x_1\\\mbf x_{2,s}\\\mbf x_{3,ss}}(t),
	\end{align*} where the 4-PI operators $\mcl{T}$, $\mcl{A}$, $\mcl{B}$, $\mcl{C}$ and $\mcl{D}$ are as defined in Eqns.~\eqref{eq:PIE_conversion}. Conversely, for given $u\in L_2([0,\infty);\R^{n_u})$ and initial conditions $\mbf v_0\in Z^{n_o,n_p}[a,b]$, suppose $\mbf v$ and $z$ satisfy the PIE defined by the 4-PI operators $\mcl{T}$, $\mcl{A}$, $\mcl{B}$, $\mcl{C}$ and $\mcl{D}$ as defined in Equations \eqref{eq:PIE_conversion}. Then, $z$ also satisfies the ODE-PDE defined by $\{A,A_i,B_i,B,B_x,C_{10},C_a,C_b,D,E,E_{10},E_a,E_b\}$ with
	\begin{align*}
		\mbf x_0 = \mcl{T}\mbf v_0,\quad \bmat{x\\\mbf x}(t) := \mcl{T}\mbf v(t).
	\end{align*}
\end{thm}
\begin{proof}
	Refer Lemma 3.3 and 3.4 in \cite{shivakumar2019generalized} for proof.
\end{proof}

PIE representations differ from typical ODE-PDE form in several significant ways. First, while PDEs rely on a differential operator in $\mcl A_d$, the a PIE system is parameterized by PI operators which are bounded on $L_2$ and form an algebra. Second, the PIE eliminates boundary conditions by incorporating the effect of boundary conditions directly into the dynamics. Finally, solutions of the PIE system are defined on $Z^{n_0,n_p}[a,b]$, which is a Hilbert space with respect to $Z$-inner product, whereas $D(\mcl{A}_d)$ is not a Hilbert space.

\section{Stability Equivalence of ODE-PDEs and PIEs}\label{sec:PIE_stability}
In this section, we show that stability in $Z$ and $L_2$-gain of the ODE-PDE system is implied by that of the PIE system. First, we define asymptotic stability of PIEs and of ODE-PDEs.
\begin{defn}
	For $w=u=0$, the PIE \eqref{eq:PIE_full} defined by $\{\mcl{T},\mcl{A},\mcl{B}_i,\mcl{C},\mcl{D}_{ij}\}$ is said to be asymptotically stable if for any initial condition $\mbf x_0\in Z^{m,n}[a,b]$, if $\mbf x$ and $z$ satisfy the PIE, we have $\lim_{t\rightarrow \infty} \norm{\mcl{T}\mbf{x}(t)}_Z= 0$.
\end{defn}

\begin{defn}
	For $w=u=0$, the ODE-PDE \eqref{eq:pde_general} defined by $\{A,A_i,B_i,B,B_x,C_{10},C_a,C_b,D,E,E_{10},E_a,E_b\}$  is said to be asymptotically stable if for initial  condition $\mbf x_0\in D(\mcl{A}_d)$, if $x$, $\mbf x$ and $z$ satisfy the ODE-PDE , then \begin{align*}
		\lim_{t\to\infty}\norm{\bmat{x\\\mbf x}(t)}_Z= 0.
	\end{align*}
\end{defn}

\begin{lem}
Suppose $\{\mcl{T},\mcl{A},\mcl{B}_i,\mcl{C},\mcl{D}_{ij}\}$ and $\{A,A_i,B_i,B,B_x,C_{10},C_a,C_b,D,E,E_{10},E_a,E_b\}$  satisfy Eqns.~\eqref{eq:PIE_conversion}. Then the ODE-PDE defined by $\{A,A_i,B_i,B,B_x,C_{10},C_a,C_b,D,E,E_{10},E_a,E_b\}$ is asymptotically stable if the PIE defined by $\{\mcl{T},\mcl{A},\mcl{B}_i,\mcl{C},\mcl{D}_{ij}\}$ is asymptotically stable.
\end{lem}
\begin{proof}
From Theorem \ref{thm:pde2pie}, $x$ and $\mbf x$ satisfy the ODE-PDE for the given $x_0$, $\mbf x_0$ if and only if $\mbf v$ satisfies the PIE for $\mbf v_0$ where
\[
\bmat{x\\\mbf x}(t)=\mcl{T}\mbf v(t),\qquad \bmat{x_0\\\mbf x_0} = \mcl T \mbf v_0.
\] 
If the PIE is stable, then $\lim_{t \rightarrow 0}\norm{\mbf{v}(t)}_{Z}=0$. Since $\mcl T$ is a bounded linear operator, this implies  
\[
\lim_{t \rightarrow 0}\norm{\bmat{x\\\mbf x}(t)}_{Z}=\lim_{t \rightarrow 0}\norm{\mcl T\mbf{v}(t)}_{Z}=0.
\]  
\end{proof}

\begin{lem}
Suppose $\{\mcl{T},\mcl{A},\mcl{B}_i,\mcl{C},\mcl{D}_{ij}\}$ and $\{A,A_i,B_i,B,B_x,C_{10},C_a,C_b,D,E,E_{10},E_a,E_b\}$ satisfy Eqns.~\eqref{eq:PIE_conversion}. For $w\in L_2([0,\infty))$, $u=0$ and $\mbf{x}(0)=0$,
	any solution $\mbf x,z$ of the PIE system satisfies $\norm{z}_{L_2}\le\gamma\norm{w}_{L_2}$ if and only if any solution to the ODE-PDE system, $x, \mbf x, z$ satisfies $\norm{z}_{L_2}\le \gamma\norm{w}_{L_2}$ for $\mbf x(0)=0$, $x(t)=0$, $u=0$ and $\bar{w}\in L_2([0,\infty))$.
\end{lem}
\begin{proof}
From Theorem \ref{thm:pde2pie}, $x$, $\mbf x$, and $z$ satisfy the ODE-PDE for the given $w$ if and only if $\mbf v$ and $z$ satisfy the PIE for the given $w$ where
\[
\bmat{x\\\mbf x}(t)=\mcl{T}\mbf v(t).
\]
\end{proof}
\section{The Dual PIE}
For a PIE system of the form Eq.~\eqref{eq:PIE_general} we may associate the following dual (adjoint) PIE.
\begin{align}
\mcl T^* \dot{\bar{\mbf x}}(t)&=\mcl A^*\bar{\mbf x}(t)+\mcl{C}^*\bar{w}(t)\notag\\
\bar{z}(t) &= \mcl{B}^*\bar{\mbf x}(t) + \mcl{D}^*\bar{w}(t)\notag
\end{align}
where $\mcl{T}^*, \mcl{A}^*:Z^{m,n}[a,b]\to Z^{m,n}[a,b]$, $\mcl{B}^*:Z^{m,n}[a,b]\to \R^{n_w}$, $\mcl{C}^*:\R^{n_z}\to Z^{m,n}[a,b]$ and $\mcl{D}^*\in \R^{n_w\times n_z}$ are 4-PI operators. 



When the PIE system Eq.\eqref{eq:PIE_general} is constructed from a PDE system, then the dual PIE system Eq.\eqref{eq:PIE_adjoint_general} may also be constructed from a PDE system. An illustrative example is given here.
\begin{ex}
	Consider the transport equation
	\begin{align}\label{eq:example_transport}
	&\dot{\mbf v}(s,t) + \mbf v_s(s,t) = 0, \qquad s\in[0,1], t>0,\notag\\
	&\mbf v(0,t) = 0, ~\mbf v(s,0) \in L_2[0,1].
	\end{align}
	The PIE form Eq.\eqref{eq:example_transport} is
	\begin{align*}
	(\threepiFull{0}{1}{0}\dot{\mbf x})(t) = (\threepiFull{-1}{0}{0}\mbf{x})(t), \qquad t>0.
	\end{align*}
	The corresponding dual PIE is
	\begin{align*}
	(\threepiFull{0}{0}{1}\dot{\mbf y})(t) = (\threepiFull{-1}{0}{0}\mbf y)(t), \qquad t>0.
	\end{align*}
	The dual PIE may be constructed from the following PDE
	\begin{align*}
	&\dot{\mbf z}(s,t) - \mbf z_s(s,t) = 0, \qquad s\in[0,1], t>0,\notag\\
	&\mbf z(1,t) = 0, ~\mbf z(s,0) \in L_2[0,1].
	\end{align*}
\end{ex}


\section{Dual Stability Theorem} \label{sec:dual}
In this section, show that stability of the dual PIE is equivalent to that of the primal PIE. 

\begin{thm}\label{thm:dual_stable}
(Dual Stability of PIEs:)	Suppose $\mcl{T}$ and $\mcl{A}$ are 4-PI operators. Then the following statements are equivalent.
	\begin{enumerate}
		\item $\lim\limits_{t\to \infty} \mcl{T}\mbf x(t)\to 0$ for any $\mbf x$ that satisfies $\mcl{T}\dot{\mbf x}(t) = \mcl{A}\mbf x(t)$ with initial condition $\mbf x(0)\in Z^{m,n}[a,b]$.
		\item $\lim\limits_{t\to \infty} \mcl{T}^*\mbf x(t)\to 0$ for any $\mbf x$ that satisfies $\mcl{T}^*\dot{{\mbf x}}(t) = \mcl{A}^*{\mbf x}(t)$ with initial condition ${\mbf x}(0)\in Z^{m,n}[a,b]$.
	\end{enumerate}
\end{thm}
\begin{proof}
	Suppose $\mbf x$ satisfies $\mcl{T}\dot{x}(t) = \mcl{A}x(t)$ with initial condition $\mbf x(0)\in Z^{m,n}[a,b]$ and $\lim_{t\rightarrow \infty}\mcl{T}\mbf{x}(t)\to 0$. Let $\bar{\mbf x}$ satisfy $\mcl{T}^*\dot{{\mbf x}}(t) = \mcl{A}^*{\mbf x}(t)$ with initial condition $\bar{\mbf x}(0)\in Z^{m,n}[a,b]$. In the following, we use $\ip{\cdot}{\cdot}=\ip{\cdot}{\cdot}_Z$. Then for any finite $t>0$, by IBP and a variable change,

{\small
		\begin{align*}
		&\int_{0}^{t} \ip{\bar{\mbf x}(t-s)}{\mcl T\dot{\mbf x}(s)}ds\\
		&=\ip{\bar{\mbf x}(0)}{\mcl T\mbf x(t)}-\ip{\bar{\mbf x}(t)}{\mcl{T}\mbf x(0)} - \int_{0}^t \ip{\partial_s{\bar{\mbf x}}(t-s)}{\mcl{T}\mbf x(s)}ds\\
		&=\ip{\bar{\mbf x}(0)}{\mcl{T}\mbf x(t)}-\ip{\bar{\mbf x}(t)}{\mcl{T}\mbf x(0)} - \int_{t}^0 \ip{\dot{\bar{\mbf x}}(\theta)}{\mcl{T}\mbf x(t-\theta)}d\theta\\
		&=\ip{\bar{\mbf x}(0)}{\mcl{T}\mbf x(t)}-\ip{\bar{\mbf x}(t)}{\mcl{T}\mbf x(0)} + \int_{0}^t \ip{\mcl T^*\dot{\bar{\mbf x}}(\theta)}{\mbf x(t-\theta)}d\theta
		\end{align*}}
	where $\theta = t-s$.
	Furthermore, using a variable change,
	\begin{align*}
	&\int_{0}^{t} \ip{\bar{\mbf x}(t-s)}{\mcl T\dot{\mbf x}(s)}ds=\int_{0}^{t} \ip{\bar{\mbf x}(t-s)}{\mcl A\mbf x(s)}ds\\
	&=\int_{0}^{t} \ip{\bar{\mbf x}(\theta)}{\mcl A\mbf x(t-\theta)}d\theta =\int_{0}^{t} \ip{\mcl A^*\bar{\mbf x}(\theta)}{\mbf x(t-\theta)}d\theta.\vspace{1mm}
	\end{align*}
	Therefore,
{\small
		\begin{align*}
		&\int_{0}^{t} \ip{\mcl A^*\bar{\mbf x}(\theta)}{\mbf x(t-\theta)}d\theta \\
		&=\ip{\bar{\mbf x}(0)}{\mcl{T}\mbf x(t)}-\ip{\bar{\mbf x}(t)}{\mcl{T}\mbf x(0)} + \int_{0}^t \ip{\mcl T^*\dot{\bar{\mbf x}}(\theta)}{\mbf x(t-\theta)}d\theta.
		\end{align*}}
	However, $\mcl A^*\bar{\mbf x}(\theta) = \mcl T^*\dot{\bar{\mbf x}}(\theta)$ for all $\theta\in [0,t]$ and so we have
	\begin{align*}
	\ip{\bar{\mbf x}(0)}{\mcl{T}\mbf x(t)}-\ip{\bar{\mbf x}(t)}{\mcl{T}\mbf x(0)}=0\qquad \forall t>0.
	\end{align*}
	Since $\lim_{t\rightarrow \infty}\mcl{T}\mbf x(t)= 0$, we have $\lim_{t\rightarrow \infty}\ip{\mcl T^*\bar{\mbf x}(t)}{\mbf x(0)}=\lim_{t\rightarrow \infty}\ip{\bar{\mbf x}(t)}{\mcl{T}\mbf x(0)}= 0$ for any $\mbf x(0)\in Z^{m,n}[a,b]$. We conclude that $\lim_{t\rightarrow \infty}\mcl T^*\bar{\mbf x}(t)= 0$. Since the dual and primal systems are interchangeable, necessity follows from sufficiency.
\end{proof}

\section{Dual $L_2$-gain Theorem}
As seen in the previous section, stability of a PIE system and its dual are equivalent. In this section, we show that for $\mbf x_0=0$, input-output performance of primal and dual PIE in the $L_2$-gain metric is equivalent.

\begin{thm}\label{thm:IO_PIE}
(Duality on $L_2$-gain bound of PIEs:)	Suppose $\mcl{T}$, $\mcl{A}$, $\mcl{B}$, $\mcl{C}$ and $\mcl{D}$ are 4-PI operators. Then the following statements are equivalent.
	\begin{enumerate}
		\item For any $w\in L_2([0,\infty);\R^q)$ and $\mbf x(0)=0$ any solution $\mbf x(t)\in Z^{m,n}$ and $z(t)\in \R^p$ of the PIE system
		\begin{align}\label{eq:PIE_IO}
		\mcl{T}\dot{\mbf x}(t) &= \mcl{A}\mbf x(t) + \mcl{B}w(t), \qquad \mbf x(0)=0\notag\\
		z(t) &= \mcl{C}\mbf x(t) + \mcl{D}w(t)
		\end{align}
		satisfies $\norm{z}_{L_2}\le \gamma\norm{w}_{L_2}$.
		\item For any $\bar{w}\in L_2([0,\infty);\R^p)$ and $\bar{\mbf x}(0)=0$, any $\bar{\mbf x}(t)\in Z^{m,n}$ and $\bar{z}(t)\in\R^q$ of the dual PIE system
		\begin{align}\label{eq:PIE_IO_adjoint}
		\mcl{T}^*\dot{\bar{\mbf x}}(t) &= \mcl{A}^*\bar{\mbf x}(t) + \mcl{C}^*\bar w(t), \qquad \bar{\mbf x}(0)=0\notag\\
		\bar{z}(t) &= \mcl{B}^*\bar{\mbf x}(t)+\mcl{D}^*\bar w(t)
		\end{align}
		satisfies $\norm{\bar{z}}_{L_2}\le \gamma\norm{\bar w}_{L_2}$.
	\end{enumerate}
\end{thm}
\begin{proof}
	Suppose that for any $w\in L_2([0,\infty);\R^q)$ and $\mbf x(0)=0$ any solution $\mbf x(t)\in Z^{m,n}$ and $z(t)\in \R^p$ of the PIE system satisfies $\norm{z}_{L_2}\le \gamma\norm{w}_{L_2}$. For $\bar{w}\in L_2([0,\infty);\R^p)$ and $\bar{\mbf x}(0)=0$, let $\bar{\mbf x}(t)\in Z^{m,n}$ and $\bar{z}(t)\in\R^q$ satisfy the dual PIE system. Then for any finite $t\ge 0$, since $\mbf x(0) = \bar{\mbf x}(0)=0$, we have 

{\small
		\begin{align*}
		&\int_{0}^{t} \ip{\bar{\mbf x}(t-s)}{\mcl T\dot{\mbf x}(s)}ds\\
		&=\ip{\bar{\mbf x}(0)}{\mcl{T}\mbf x(t)}-\ip{\bar{\mbf x}(t)}{\mcl{T}\mbf x(0)} + \int_{0}^t \ip{\mcl T^*\dot{\bar{\mbf x}}(\theta)}{\mbf x(t-\theta)}d\theta\\
		& = \int_{0}^t \ip{\mcl T^*\dot{\bar{\mbf x}}(\theta)}{\mbf x(t-\theta)}d\theta
		\end{align*}}
	where $\theta = t-s$. Furthermore, by the change variable change,
	\begin{align*}
	&\int_{0}^{t} \ip{\bar{\mbf x}(t-s)}{\mcl T\dot{\mbf x}(s)}ds\\
	&=\int_{0}^{t} \ip{\bar{\mbf x}(t-s)}{\mcl A\mbf x(s)}ds+\int_{0}^{t} \ip{\bar{\mbf x}(t-s)}{\mcl Bw(s)}ds\\
	&=\int_{0}^{t} \ip{\mcl A^*\bar{\mbf x}(\theta)}{\mbf x(t-\theta)}d\theta+\int_{0}^{t} \ip{\mcl B^*\bar{\mbf x}(\theta)}{w(t-\theta)}d\theta.
	\end{align*}
	Combining, we obtain
	\begin{align*}
	&\int_{0}^t \ip{\mcl T^*\dot{\bar{\mbf x}}(\theta)}{\mbf x(t-\theta)}d\theta\\
	& = \int_{0}^{t} \ip{\mcl A^*\bar{\mbf x}(\theta)}{\mbf x(t-\theta)}d\theta+ \int_{0}^{t} \ip{\mcl B^*\bar{\mbf x}(\theta)}{w(t-\theta)}d\theta.		
	\end{align*}
Now, by the definition of $\bar z$, we obtain
\begin{align*}
	&\int_{0}^{t} \ip{\bar{z}(\theta)}{ w(t-\theta)}d\theta-\int_{0}^{t} \ip{\mcl D^*\bar{w}(\theta)}{w(t-\theta)}d\theta\\
	&=\int_{0}^{t} \ip{\mcl B^*\bar{\mbf x}(\theta)}{ w(t-\theta)}d\theta\\
	&=\int_{0}^t \ip{\mcl T^*\dot{\bar{\mbf x}}(\theta)}{\mbf x(t-\theta)}d\theta- \int_{0}^{t} \ip{\mcl A^*\bar{\mbf x}(\theta)}{\mbf x(t-\theta)}d\theta\\
	&=\int_{0}^t \ip{\mcl C^*\bar{w}(\theta)}{\mbf x(t-\theta)}d\theta=\int_{0}^t \ip{\bar{w}(\theta)}{\mcl C\mbf x(t-\theta)}d\theta\\
	&=\int_{0}^t \ip{\bar{w}(\theta)}{z(t-\theta)}d\theta-\int_{0}^t \ip{\bar{w}(\theta)}{\mcl Dw(t-\theta)}d\theta.
	\end{align*}
	We conclude that for any $t>0$, if $z$ and $w$ satisfy the primal PIE and  $\bar{z}$ and $\bar w$ satisfy the dual PIE, then
	\begin{align*}
	\int_{0}^{t} \ip{\bar{z}(\theta)}{ w(t-\theta)}d\theta = \int_{0}^t \ip{\bar{w}(\theta)}{z(t-\theta)}d\theta.
	\end{align*}
Now, for any $\bar w \in L_2^p$, suppose $\bar z$ solves the dual PIE for some $\bar{\mbf x}$. For any fixed $T>0$, define $w(t) = \bar{z}(T-t)$ for $t\le T$ and $w(t)=0$ for $t >T$. Then $w \in L_2^q$ and for this input, let $z$ solve the primal PIE for some $\mbf x$. Then if we define the truncation operator $P_T$, we have
	\begin{align*}
	\norm{P_T\bar{z}}_{L_2}^2 &= \int_{0}^T \ip{\bar{z}(s)}{\bar{z}(s)}ds= \int_{0}^T \ip{\bar{z}(s)}{w(T-s)}ds \\
	&= \int_{0}^T \ip{\bar{w}(s)}{z(T-s)}ds\le \norm{P_T\bar{w}}_{L_2}\norm{P_Tz}_{L_2}\\
	&\le\gamma \norm{\bar{w}}_{L_2}\norm{w}_{L_2} = \gamma\norm{\bar{w}}_{L_2}\norm{\bar{z}}_{L_2}.
	\end{align*}
	\begin{align*}
	\norm{P_T\bar{z}}_{L_2}^2 &= \int_{0}^T \ip{\bar{z}(s)}{\bar{z}(s)}ds= \int_{0}^T \ip{\bar{z}(s)}{w(T-s)}ds \\
	&= \int_{0}^T \ip{\bar{w}(s)}{z(T-s)}ds\le \norm{P_T\bar{w}}_{L_2}\norm{P_Tz}_{L_2}\\
	&\le\norm{P_T\bar{w}}_{L_2}\norm{z}_{L_2}\le\gamma \norm{P_T\bar{w}}_{L_2}\norm{w}_{L_2}\\
&= \gamma\norm{P_T\bar{w}}_{L_2}\norm{P_T w}_{L_2} = \gamma\norm{P_T\bar{w}}_{L_2}\norm{P_T\bar{z}}_{L_2}.
	\end{align*}
Therefore, we have that $\norm{P_T\bar{z}}_{L_2} \le \gamma  \gamma\norm{P_T\bar{w}}_{L_2}$ for all $T\ge 0$. Hence, we conclude that $\norm{\bar{z}}_{L_2} \le \gamma  \gamma\norm{\bar{w}}_{L_2}$. Since the dual and primal systems are interchangeable, necessity follows from sufficiency.
\end{proof}
\section{Linear Partial Integral Inequalities}\label{sec:LPIs}
Optimization problems with PI operator decision variables and Linear PI Inequality constraints are called Linear PI Inequalities (LPIs) and take the form
\begin{align}\label{eq:LPI}
&\fourpi{P_0}{Q_0}{Q_0^T}{R_{0i}}+\sum_{k=1}^Nx_j\fourpi{P_k}{Q_k}{Q_k^T}{R_{ki}}\succcurlyeq 0,
\end{align}
where the decision variable is $x\in\R^N$ and $\fourpi{P_k}{Q_k}{Q_k^T}{R_{ki}}:\R^m\times L_2^n[a,b]\to \R^m\times L_2^n[a,b]$ is a given self-adjoint 4-PI operator for $0\le k\le N$ and $k\in\Z$.

LPI optimization problems can be solved using the MATLAB software package PIETOOLS~\cite{toolbox:pietools}. In the following sections, we present applications of Theorems \ref{thm:dual_stable} and \ref{thm:IO_PIE} in the form of LPI tests for dual stability, dual $L_2$-gain, stabilization, and $\hinf$-optimal control of PIE systems, each with associated code snippets using the PIETOOLS implementation.
\section{A Dual LPI For Stability}

Using Theorem \ref{thm:dual_stable}, we give primal and dual LPIs for stability of a PIE system.

\begin{thm}\label{thm:primal_LPI}
	(Primal LPI for Stability:)	Suppose there exists a self-adjoint bounded and coercive operator $\mcl{P}:Z^{m,n}[a,b] \to Z^{m,n}[a,b]$ such that \begin{align}
	\mcl{T}^*\mcl{P}\mcl{A}+\mcl{A}^*\mcl{P}\mcl{T}\preccurlyeq -\epsilon \mcl{T}^*\mcl{T} \label{eq:primal_stable_LPI}
	\end{align} for some $\epsilon>0$. Then any $\mbf{x}\in Z^{m,n}[a,b]$ that satisfies the system \begin{align*}
	\mcl{T}\dot{\mbf x}(t) = \mcl{A}\mbf x(t), \quad \mbf{x}(0)=\mbf x_0\in Z^{m,n}[a,b]
	\end{align*} we have $\norm{\mcl{T}\mbf x(t)}_{Z}\le\norm{\mcl{T}\mbf x(0)}_{Z} M e^{(-\alpha t)}$ for some $M$ and $\alpha>0$.
\end{thm}
\begin{proof}
	The proof can be found in the \cite{peet_2020Aut}.
\end{proof}

\begin{thm}\label{thm:dual_LPI}
(Dual LPI for Stability:)	Suppose there exists a self-adjoint bounded and coercive operator $\mcl{P}:Z^{m,n}[a,b] \to Z^{m,n}[a,b]$ such that \begin{align}
 \mcl{T}\mcl{P}\mcl{A}^*+\mcl{A}\mcl{P}\mcl{T}^*\preccurlyeq -\epsilon \mcl{T}\mcl{T}^* \label{eq:dual_stable_LPI}
\end{align} for some $\epsilon>0$. Then any $\mbf{x}\in Z^{m,n}[a,b]$ that satisfies the system \begin{align*}
	\mcl{T}\dot{\mbf x}(t) = \mcl{A}\mbf x(t), \quad \mbf{x}(0)=\mbf x_0\in Z^{m,n}[a,b]
\end{align*} we have $\lim_{t\rightarrow \infty}\norm{\mcl{T}\mbf x(t)}_{Z}= 0$.
\end{thm}
\begin{proof}
	The proof can be found in the Appendix.
\end{proof}
{\small
\begin{pse}
	\hfill
	\begin{flalign*}
	&\texttt{prog = sosprogram([s,t]);}&\\
	&\texttt{[prog,P] = sos\_posopvar(prog,dim,I,s,t);}&\\ %
	&\texttt{D=T*P*A'+A*P*T'+eps*T*T';}&\\
	&\texttt{prog = sos\_opineq(prog, -D);}&\\
	&\texttt{prog = sossolve(prog);}&
	\end{flalign*}
\end{pse}
}
\section{Dual KYP Lemma}
We formulate the following dual LPI for $L_2$-gain of PIE system in the form Eq.~\eqref{eq:PIE_general} where $\mcl{T}:Z^{m,n}[a,b] \to Z^{m,n}[a,b]$, $\mcl{A}:Z^{m,n}[a,b] \to Z^{m,n}[a,b]$, $\mcl{B}_i:\R^q\to Z^{m,n}[a,b]$, $\mcl{C}:Z^{m,n}[a,b] \to \R^r$ and $\mcl{D}_{1i}:\R^q\to\R^r$.
\begin{thm}\label{thm:gain_LPI}
	(LPI for $L_2$-gain:) Suppose there exist $\epsilon>0, \gamma>0$, bounded linear operators $\mcl{P}:Z^{m,n}[a,b] \to Z^{m,n}[a,b]$, such that $\mcl{P}$ is self-adjoint, coercive and
	\begin{align}\label{eq:LPI_gain}
	&\bmat{-\gamma I&\mcl D&\mcl{C}\mcl{P}\mcl{T}^*\\
		(\cdot)^*&-\gamma I&\mcl B^*\\
		(\cdot)^*&(\cdot)^*&(\cdot)^*+\mcl T(\mcl{A}\mcl{P})^*}	\preccurlyeq 0.
	\end{align}
	Then, for $w\in L_2$, any $\mbf x$ and $z$ that satisfy the PIE \eqref{eq:PIE_general} also satisfies $\norm{z}_{L_2}\le \gamma \norm{w}_{L_2}$.
\end{thm}
\begin{proof}
	The proof is same as the proof for Theorem \ref{thm:IO_LPI} with $\mcl{B}_1=\mcl{B}$, $\mcl{B}_2=0$, $\mcl{D}_{11}=\mcl{D}$ and $\mcl{D}_{12}=0$.
\end{proof}
{\small
\begin{pse}
	\hfill
\begin{flalign*}
&\texttt{prog = sosprogram([s,t], gam);}&\\
&\texttt{[prog,P] = sos\_posopvar(prog,dim,I,s,t);}&\\ %
&\texttt{D = [-gam*I+eps*I \hspace{1mm} D \hspace{5.5mm} C*P*T';}&\\
&\texttt{\hspace{9mm} D \hspace{8mm} -gam*I+eps*I \hspace{0mm} B';}&\\
&\texttt{\hspace{7mm} (P*C')' \hspace{0mm} B' \hspace{0mm} ($\cdot$)'+T*(A*P)'+eps*T*T'];}&\\
&\texttt{prog = sos\_opineq(prog, -D);}&\\
&\texttt{prog = sossetobj(prog,gam);}&\\
&\texttt{prog = sossolve(prog);}&
\end{flalign*}
\end{pse}
}
\section{Stabilizing Controller Synthesis}
For PIEs with inputs,
\[\mcl{T}\dot{\mbf x}(t) = \mcl{A}\mbf x(t)+\mcl{B}u(t)\] the following LPI can be used to find a stabilizing state-feedback controller of the form $u(t) = \mcl{K}\mbf{x}(t)$ where $\mcl{K}:Z^{m,n}[a,b]\to \R^q$ is a 4-PI operator.
\begin{cor}\label{thm:dual_stabcon_LPI}
(LPI for Stabilizing Controller Synthesis:)	Suppose there exist bounded linear operators $\mcl{P}:Z^{m,n}[a,b] \to Z^{m,n}[a,b]$ and $\mcl{Z}:Z^{m,n}[a,b]\to \R$, such that $\mcl{P}$ is self-adjoint, coercive and
	\begin{align}
	(\mcl{AP+BZ})\mcl{T}^*+\mcl{T}(\mcl{AP+BZ})^*\le -\epsilon \mcl{T}\mcl{T}^*. \label{eq:dual_control}\end{align}
	Then, for $u(t) = \mcl{K}\mbf x(t)$, where $\mcl{K}=\mcl{ZP}^{-1}$, any $\mbf{x}\in Z^{m,n}[a,b]$ that satisfies the system \[\mcl{T}\dot{\mbf x}(t) = \mcl{A}\mbf x(t)+\mcl{B}u(t), \quad \mbf{x}(0)=\mbf x_0\in Z^{m,n}[a,b]\] also satisfies $\lim_{t\rightarrow \infty}\norm{\mcl{T}\mbf x(t)}_{Z}= 0$.
\end{cor}
\begin{proof}
	The proof is same as the proof for Theorem \ref{thm:dual_LPI} substituting $\mcl{A}\rightarrow \mcl{A}+\mcl{BK}$ and where $\mcl{Z}=\mcl{KP}$.
\end{proof}
{\small
\begin{pse}
	\hfill
	\begin{flalign*}
	&\texttt{prog = sosprogram([s,t]);}&\\
	&\texttt{[prog,Z] = sos\_opvar(prog,dim,I,s,t,deg);}&\\
	&\texttt{[prog,P] = sos\_posopvar(prog,dim,I,s,t);}&\\ %
	&\texttt{D=T*(A*P+B*Z)'+(A*P+B*Z)*T'+eps*T*T';}&\\
	&\texttt{prog = sos\_opineq(prog, -D);}&\\
	&\texttt{prog = sossolve(prog);}&
	\end{flalign*}
\end{pse}
}
\section{$\hinf$-optimal Controller Synthesis}
For PIE systems with inputs and outputs, we can use Theorem \ref{thm:IO_PIE} to pose the $\hinf$-optimal controller synthesis problem as an LPI. Specifically, we formulate the following LPI for finding the $\hinf$-optimal controller for a PIE system in the form Eq.~\eqref{eq:PIE_full} where $\mcl{T}, \mcl{A}:Z^{m,n}[a,b] \to Z^{m,n}[a,b]$, $\mcl{B}_1:\R^q\to Z^{m,n}[a,b]$, $\mcl{B}_2:\R^p\to Z^{m,n}[a,b]$, $\mcl{C}:Z^{m,n}[a,b]\to \R^r$, $\mcl{D}_{11}:\R^q\to\R^r$ and $\mcl{D}_{12}:\R^p\to\R^r$.
\begin{thm}\label{thm:IO_LPI}
(LPI for $H_{\infty}$ Optimal Controller Synthesis:) Suppose there exist $\gamma>0$, bounded linear operators $\mcl{P}:Z^{m,n}[a,b]\to Z^{m,n}[a,b]$ and $\mcl{Z}:Z^{m,n}[a,b]\to \R^p$, such that $\mcl{P}$ is self-adjoint, coercive and
	\begin{align}\label{eq:LPI_con}
	&\bmat{-\gamma I&\mcl D_{11}&(\mcl{C}\mcl{P}+\mcl{D}_{12}\mcl{Z})\mcl{T}^*\\
		(\cdot)^*&-\gamma I&\mcl B_1^*\\
		(\cdot)^*&(\cdot)^*&(\cdot)^*+\mcl T(\mcl{A}\mcl{P}+\mcl{B}_2\mcl{Z})^*}	\preccurlyeq 0.
	\end{align}
	Then, for any $w\in L_2$, for $u(t) = \mcl{K}\mbf x(t)$ where $\mcl{K}=\mcl{ZP}^{-1}$, any $\mbf x$ and $z$ that satisfy the PIE \eqref{eq:PIE_full} also satisfy $\norm{z}_{L_2}\le \gamma \norm{w}_{L_2}$.
\end{thm}
\begin{proof}
	The proof can be found in the appendix.
\end{proof}

{\small
\begin{pse}
	\hfill
\begin{flalign*}
&\texttt{prog = sosprogram([s,t], gam);}&\\
&\texttt{[prog,P] = sos\_posopvar(prog,dim,I,s,t);}&\\ %
&\texttt{[prog,Z] = sos\_opvar(prog,dim,I,s,t,deg);}&\\
&\texttt{D = [-gam*I+eps*I \hspace{1mm} D11' \hspace{5.5mm} (C*P+D12*Z)*T';}&\\
&\texttt{\hspace{6mm} ($\cdot$)' \hspace{8mm} -gam*I+eps*I \hspace{0mm} B1';}&\\
&\texttt{\hspace{6mm} ($\cdot$)' \hspace{0mm} ($\cdot$)' \hspace{0mm} ($\cdot$)'+T*(A*P + B2*Z)'+eps*T*T'];}&\\
&\texttt{prog = sos\_opineq(prog, -D);}&\\
&\texttt{prog = sossetobj(prog,gam);}&\\
&\texttt{prog = sossolve(prog);}&
\end{flalign*}
\end{pse}
}
\section{Numerical Examples}

In this section, use various numerical examples to demonstrate the accuracy and scalability of the LPIs presented in this paper. First, we verify the stability of PDEs, where the stability holds for certain values of the system parameters (referred to as a stability parameter). We test for the stability of the system using the dual stability criterion and change the stability parameter continuously to identify the point at which the stability of the system changes. The second set of examples will focus on finding in-domain controllers to stabilize an unstable system. Finally, we also present a numerical example of systems with inputs and outputs to find $\hinf$-optimal controllers.

\subsection{Stability Tests Using Dual Stability Criterion}
\begin{ex}
	Consider the scalar diffusion-reaction equation with fixed boundary conditions.
	\begin{align*}
	&u_t(s,t) = \lambda u(s,t) + u_{ss}(s,t), \qquad s\in[0,1], t>0,\\
	& u(0,t)=u(1,t)=0, u(s,0) = u_0
	\end{align*}
	We can establish analytically that this system is stable for $\lambda \le \pi^2$. We increase $\lambda$ continuously and determine the maximum value for which the system is stable. From our tests, we find that the system is stable for $\lambda \le (1+1e^{-5})\pi^2$.
\end{ex}
\begin{ex}
	Let us change the boundary conditions of the previous example. Then the bound on the stability parameter changes to $\lambda\le 2.467$.
	\begin{align*}
	&u_t(s,t) = \lambda u(s,t) + u_{ss}(s,t), \qquad s\in[0,1], t>0,\\
	& u(0,t)=u_s(1,t)=0, u(s,0) = u_0
	\end{align*}
	Testing for stability using the dual lyapunov criterion we find that the system is stable for $\lambda \le 2.467+5e^{-4}$.
\end{ex}

\subsection{Finding Stabilizing Controller For Unstable PDE Systems}

\begin{ex}
	Let us revisit the Example 8. Suppose $\lambda=10$. Then the system is unstable. To stabilize the system, we introduce an in-domain control input as
	\[
	u_t(s,t) = \lambda u(s,t) + u_{ss}(s,t) + d(t)
	\]
	where $d(t) = \myint K(s) u_{ss}(s,t) ds$ is the control input. Solving the LPI in Theorem \ref{thm:dual_stabcon_LPI} we get the controller
	\begin{align*}
		K(s) =& 0.29s^5-1.01s^4+0.95s^3+0.16s^2-0.51s+0.98.
	\end{align*}
	
\end{ex}

\subsection{$\hinf$-optimal Controller Synthesis}	
\begin{ex} \label{ex:hinf}
	Consider the following cascade of diffusion-reaction equations with a dynamic controller acting at the boundary.
	{\small
		\begin{align*}
		&\dot{\mbf x}_i(s,t) = \lambda\mbf  x_i(s,t) +\sum_{k=i}^{N}\mbf  x_{k,ss}(s,t)+w(t),\quad  i\in\{1...N\}\\
		&\mbf x_i(0,t) =0, \mbf x_{i}(1,t) = 0 \quad\forall ~i\in\{1...N-1\},\\
		&\mbf x_N(0,t)=0, \mbf  x_N(1,t) = x_0(t)\\
		&\dot{x}_0(t) = u(t), \qquad x_0(0)=0, \mbf x_i(s,0)=0, s\in[0,1]\\
		&z(t) = x_0(t)
		\end{align*}}
	where $x_0$ is the state of the dynamic boundary controller, $\mbf{x}_i$ are distributed states, $z$ is the output and $w_i$ are the input disturbances. The control input, $u(t) = K_0x_0(t)+\int_0^1 K(s) \mbf x(s,t) ds$ where $K:[a,b]\to \R^{1\times N}$, enters the system through the ODE and acts at the boundary of the PDE state $x_N$.
	For $\lambda=10$, $N=3$ the $\hinf$-optimal controller has a norm bound of $6.5095$. In Figure 1, we plot the system response for a disturbance $w(t) = \frac{sin(5t)}{3t}$ with zero initial conditions.
\end{ex}

\begin{figure}
	\centering
	\includegraphics*[scale=0.5]{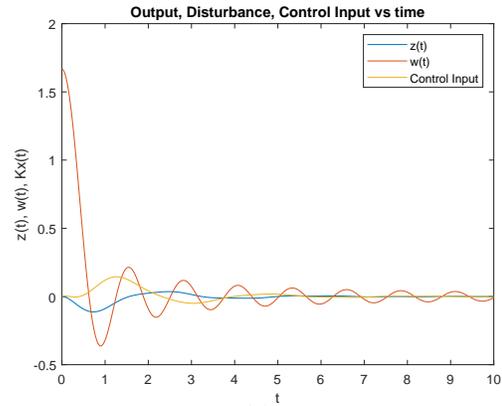}\vspace{-5mm}
	\caption{The plot shows variation of Output $z$ and control input $u$ with time when a bounded disturbance $w$ is applied to system in Example \ref{ex:hinf} with zero initial conditions.}\vspace{-5mm}
	\label{fig:plot1}
\end{figure}

\section{CONCLUSIONS}\vspace{-2mm}
In this article, we have proven the equivalence, in stability and $\hinf$-norm, between a PIE system and its dual system. Coupled ODE-PDE systems have equivalent PIE representations and properties of the ODE-PDE system are inherited from the PIE. Our duality results allow us to use Linear PI Inequalities to find stabilizing and $\hinf$-optimal state-feedback controllers for PIE systems and these controllers can then be used to regulate the associated ODE-PDE systems. We have demonstrated the accuracy and scalability of the resulting algorithms by applying the results to several illustrative examples. While the scope of the paper is limited to inputs entering through the ODE or in-domain, we believe the results can be extended to inputs at the boundary.\vspace{-2mm}

\section*{ACKNOWLEDGMENT}
This work was supported by Office of Naval Research
Award N00014-17-1-2117.\vspace{-2mm}

\bibliography{references, peet_bib}
\bibliographystyle{plain}

\section*{APPENDIX}

\subsection{Proof of Theorem \ref{thm:dual_LPI}}
\setcounter{thm}{12}
\begin{thm}
	(Dual LPI for Stability:)	Suppose there exists a self-adjoint bounded and coercive operator $\mcl{P}:Z^{m,n}[a,b] \to Z^{m,n}[a,b]$ such that \begin{align*}
	\mcl{T}\mcl{P}\mcl{A}^*+\mcl{A}\mcl{P}\mcl{T}^*\preccurlyeq -\epsilon \mcl{T}\mcl{T}^*
	\end{align*} for some $\epsilon>0$. Then any $\mbf{x}\in Z^{m,n}[a,b]$ that satisfies the system $$\mcl{T}\dot{\mbf x}(t) = \mcl{A}\mbf x(t), \quad \mbf{x}(0)=\mbf x_0\in Z^{m,n}[a,b]$$ satisfies $\norm{\mcl{T}\mbf x(t)}_{Z}\to 0$ as $t\to\infty$.
\end{thm}
\begin{proof}
	Define a Lyapunov candidate as $V(y) = \ip{\mcl{T}^*y}{\mcl{P}\mcl{T}^*y}_Z$. Then there exists an $\alpha>0$ and $\beta>0$ such that
	\begin{align*}
	\alpha\norm{\mcl{T}^*y}_Z\le V(y) \le \beta \norm{\mcl{T}^*y}_Z.
	\end{align*}
	The time derivative of $V(y)$ along the solutions of the PIE \begin{align*}\mcl{T}^*\dot{y}(t)=\mcl{A}^*y(t), \qquad y(0)\in Z^{m,n}[a,b]
	\end{align*} is given by
	\begin{align*}
	\dot{V}(y(t)) &= \ip{\mcl{T}^*y(t)}{\mcl{P}\mcl{T}^*\dot y(t)}_Z+\ip{\mcl{T}^*\dot y(t)}{\mcl{P}\mcl{T}^*y(t)}_Z\\
	&= \ip{\mcl{T}^*y(t)}{\mcl{P}\mcl{A}^*y(t)}_Z+\ip{\mcl{A}^* y(t)}{\mcl{P}\mcl{T}^*y(t)}_Z\\
	&= \ip{y(t)}{\mcl{T}\mcl{P}\mcl{A}^*y(t)}_Z+\ip{y(t)}{\mcl{A}\mcl{P}\mcl{T}^*y(t)}_Z\\
	&\le -\epsilon\norm{\mcl{T}^*y(t)}_Z\le -\frac{\epsilon}{\beta}V(y(t)).
	\end{align*}
	Then, by using Gronwall-Bellman Inequality, there exists constants $M$ and $k$ such that
	\begin{align*}
	V(y(t))\le V(y(0))Me^{(-kt)}.
	\end{align*}
	As $t\to\infty$, $V(y(t))\to0$ which implies $\norm{\mcl{T}^*y(t)}_Z\to 0$. Then, from Theorem \ref{thm:dual_stable}, $\norm{\mcl{T}x(t)}_Z\to 0$.
\end{proof}
\subsection{Proof of Theorem \ref{thm:IO_LPI}}
\setcounter{thm}{15}
\begin{thm}
	(LPI for $H_{\infty}$ Optimal Controller Synthesis:)		If there exist $\gamma>0$, bounded linear operators $\mcl{P}:Z^{m,n}[a,b]\to Z^{m,n}[a,b]$ and $\mcl{Z}:Z^{m,n}[a,b]\to \R^p$, such that $\mcl{P}$ is self-adjoint, coercive and
	\begin{align}
	&\bmat{-\gamma I&\mcl D_{11}&(\mcl{C}\mcl{P}+\mcl{D}_{12}\mcl{Z})\mcl{T}^*\\
		(\cdot)^*&-\gamma I&\mcl B_1^*\\
		(\cdot)^*&(\cdot)^*&(\cdot)^*+\mcl T(\mcl{A}\mcl{P}+\mcl{B}_2\mcl{Z})^*}\preccurlyeq 0.
	\end{align}
	then for $u = \mcl{K}\mbf x$, where $\mcl{K}=\mcl{ZP}^{-1}$, and any $w\in L_2$, any $\mbf x$ and $z$ that satisfy the PIE \eqref{eq:PIE_full} also satisfies $\norm{z}_{L_2}\le \gamma \norm{w}_{L_2}$.
\end{thm}
\begin{proof}
	Define a Lyapunov candidate function $V(x) = \ip{\mcl{T}^*x}{\mcl{P}\mcl{T}^*x}_Z$. Since $\mcl{P}$ is coercive and bounded, there exists $\alpha>0$ and $\beta>0$ such that
	\begin{align*}
	\alpha\norm{\mcl{T}^*x}_Z\le V(x) \le \beta \norm{\mcl{T}^*x}_Z.
	\end{align*}
	The time derivative of $V(x)$ along the solutions of
	\begin{align}\label{eq:generalPIE}
	\mcl{T}^*\dot{\mbf x}(t) &= (\mcl{A}+\mcl{B}_2\mcl{K})^*\mbf x(t) + (\mcl{C}+\mcl{D}_{12}\mcl{K})^*w(t), \notag\\
	z(t) &= \mcl{B}_1^*\mbf x(t) + \mcl{D}_{11}^*w(t), \quad\mbf x(0)=0
	\end{align} is given by
	\begin{align*}
	\dot{V}(x(t)) &= \ip{\mcl{T}^*x(t)}{\mcl{P}\mcl{T}^*\dot x(t)}_Z+\ip{\mcl{T}^*\dot x(t)}{\mcl{P}\mcl{T}^*x(t)}_Z\\
	&= \ip{\mcl{T}^*x(t)}{\mcl{P}(\mcl{A}+\mcl{B}_2\mcl{K})^*x(t)}_Z\\
	&\quad+\ip{(\mcl{A}+\mcl{B}_2\mcl{K})^* x(t)}{\mcl{P}\mcl{T}^*x(t)}_Z\\
	&~~\quad+ \ip{\mcl{T}^*x(t)}{\mcl{P}(\mcl{C}+\mcl{D}_{12}\mcl{K})^*w(t)}_Z\\
	&~~~~\quad+\ip{(\mcl{C}+\mcl{D}_{12}\mcl{K})^* w(t)}{\mcl{P}\mcl{T}^*x(t)}_Z.
	\end{align*}
	For any $w(t)\in\R^p$ and $x(t)\in Z$ that satisfies Eq. \eqref{eq:generalPIE},
	{\normalsize
	\begin{align*}
	&{\small\ip{\bmat{v(t)\\w(t)\\\mbf x(t)}}{\bmat{-\gamma I&\mcl D_{11}&(\mcl{C}\mcl{P}+\mcl{D}_{12}\mcl{Z})\mcl{T}^*\\
				(\cdot)^*&-\gamma I&\mcl B_1^*\\
				(\cdot)^*&(\cdot)^*&(\cdot)^*+\mcl T(\mcl{A}\mcl{P}+\mcl{B}_2\mcl{Z})^*}\bmat{v(t)\\w(t)\\\mbf x(t)}}}\\
	&= \ip{\bmat{w(t)\\\mbf x(t)}}{\bmat{0&(\mcl{C}\mcl{P}+\mcl{D}_{12}\mcl{Z})\mcl{T}^*\\(\cdot)^*&(\cdot)^*+\mcl{T}(\mcl{A}\mcl{P}+\mcl{B}_2\mcl{Z})^*}\bmat{w(t)\\\mbf x(t)}}\\
	&\quad -\gamma\norm{w(t)}^2 -\gamma\norm{v(t)}^2+\ip{v(t)}{\mcl{B}_1^*\mbf x(t)+\mcl D_{11}^* w(t)}\\
	&\quad+\ip{\mcl{B}_1^*\mbf x(t)+\mcl D_{11}^*w(t)}{v(t)}\\
	&= \dot{V}(\mbf x(t)) -\gamma\norm{w(t)}^2 -\gamma\norm{v(t)}^2+\ip{v(t)}{z(t)}\\
	&\quad+\ip{z(t)}{v(t)}\le 0
	\end{align*}}
	for any $v(t)\in\R^{p}$ and $t\ge0$. Let $v(t) = \frac{1}{\gamma}z(t)$. Then
	\begin{align*}
	\dot{V}(\mbf x(t))&\le \gamma\norm{w(t)}^2 +\frac{1}{\gamma}\norm{z(t)}^2-\frac{2}{\gamma}\norm{z(t)}^2 -\epsilon \norm{\mbf x(t)}^2\\
	&\le \gamma\norm{w(t)}^2 -\frac{1}{\gamma}\norm{z(t)}^2.
	\end{align*}
	Integrating forward in time with the initial condition $\mbf x(0)=0$, we get
	\begin{align*}
	\frac{1}{\gamma}\norm{z(t)}^2\le\gamma\norm{w(t)}^2.
	\end{align*}
	Using Theorem \ref{thm:IO_PIE}, the adjoint PIE system of Eq.\eqref{eq:generalPIE} has the same bound on $L_2$-gain from input to output. In other words, for $w\in L_2$, any $\mbf{x}$ and $z$ that satisfy equations \begin{align*}
	\mcl{T}\dot{\mbf x}(t) &= \mcl{A}\mbf x(t) + \mcl{B}_1w(t)+\mcl{B}_2u(t), \qquad \mbf x(0)=0\notag\\
	z(t) &= \mcl{C}\mbf x(t) + \mcl{D}_{11}w(t)+\mcl{D}_{12}u(t)
	\end{align*} with $u = \mcl{K}\mbf x$ and $\mcl{K}=\mcl{Z}\mcl{P}^{-1}$, we have $\norm{z}_{L_2}\le \gamma\norm{w}_{L_2}$.
\end{proof}

\subsection{PI operator definitions in Theorem \ref{thm:pde2pie}}
\newcommand{\fourpiA}[4]{\ensuremath{\mcl{P}{\scriptsize\bmat{#1,& \hspace{-3mm}#2 \\ #3,& \hspace{-3mm} \left\{#4\right\} }}}}
\newcommand{\fourpiFullA}[6]{\mcl{P}{\scriptsize\bmat{#1,& \hspace{-3mm}#2 \\ #3,& \hspace{-3mm} \left\{#4,#5,#6\right\} }}}
Given $\{A,A_i,B_i,B,B_x,C_{10},C_a,C_b,D,E,E_{10},E_a,E_b\}$, we define the following functions and 4-PI operators.

{\scriptsize\begin{align}
&H_0(s) = K(s)(BT)^{-1}B_x, H_1(s) = V(s)(BT)^{-1}B_x,\notag\\
&T_1 = T(BT)^{-1}B_x, \quad T_2(s) = T(BT)^{-1}BQ(0,s)+Q(0,s)\notag\\
&G_0 = \bmat{I & 0 & 0\\0 &0&0\\0 &0&0},~ L_0 = \bmat{0 & I & 0\\0 &0&0},\notag\\
&G_2(s,\theta) = -K(s)(BT)^{-1}BQ(s,\theta),\notag\\
&G_1(s,\theta) = \bmat{0 & 0 & 0\\0 &I&0\\0 &0&(s-\theta)I}+G_2(s,\theta),\notag\\
&L_2(s,\theta) = -V(s)(BT)^{-1}BQ(s,\theta),\; L_1(s,\theta) = \bmat{0 & 0 & 0\\0 &0&I}+L_2(s,\theta),\notag\\
&K(s) = \bmat{0 & 0 & 0\\I &0&0\\0 &I&(s-a)I},\quad V(s) = \bmat{0 & 0 & 0\\0 &0&I}\notag\\
&T = \bmat{I &0&0\\0 & I & 0\\0&0&I\\I &0&0\\0 &I&(b-a)I\\0&0&I}, Q(s,\theta) = \bmat{0 &0&0\\0 & 0 & 0\\0&0&0\\0 &I&0\\0 &0&(b-\theta)I\\0&0&I}\notag\\
&\mcl{T} = \fourpiA{I}{0}{H_0}{G_i}, \mcl{B} = \fourpiA{B_1}{0}{B_2}{0}, \mcl{D} = \fourpiA{D}{0}{0}{0}\notag\\
&\mcl{A} = \fourpiFullA{A+E_{10}T_1}{E_{10}T_2}{E}{\bmat{0~0~A_2}}{0}{0} \notag\\
&+ \fourpiFullA{0}{E_a}{0}{A_0}{0}{0}\fourpiA{0}{0}{H_0}{G_i}\notag + \fourpiFullA{0}{E_b}{0}{A_1}{0}{0}\fourpiA{0}{0}{H_1}{L_i}\notag\\
&\mcl{C} = \fourpiA{C+C_{10}T_1}{C_{10}T_2}{0}{0}+\fourpiA{0}{C_a}{0}{0}\fourpiA{0}{0}{H_0}{G_i}\notag\\
&\qquad+\fourpiA{0}{C_b}{0}{0}\fourpiA{0}{0}{H_1}{L_i}. \label{eq:PIE_conversion}
\end{align}}

\vfill
\end{document}